\documentclass[10pt,leqno]{amsart}
\usepackage{graphicx}
\usepackage{indentfirst,csquotes}

\usepackage{amssymb,amsthm,amsmath}
\usepackage{xcolor,paralist,hyperref,titlesec,fancyhdr,etoolbox}

\titleformat{\section}[display]{\normalfont\huge\bfseries\centering}{\centering\chaptertitlename\thechapter}{10pt}{\Large}
\titlespacing*{\section}{0pt}{0ex}{0ex}

\hypersetup{ colorlinks=true, linkcolor=black, filecolor=black, urlcolor=black }

\usepackage{lipsum}

\begin{document}
\title{Microscopic limits of PDEs modeling macroscopic heat conduction}

%\author{
%\name{Ronald Mickens\textsuperscript{a} and Talitha Washington\textsuperscript{b}\thanks{CONTACT T.~W. Author. Email: twashington@cau.edu}}
%\affil{\textsuperscript{a}Department of Physics,  Clark Atlanta University,  Atlanta,  GA, USA;\textsuperscript{b}Department of Mathematical Sciences,  Clark Atlanta University,  Atlanta,  GA,  USA}}
%\author[Initial Surname]{Author}
%\date{\today}
%\address{Address}
%\email{example@mail.com}

\author[Mickens]{Ronald Mickens}
%\date{\today}
\address{Department of Physics, Clark Atlanta University, Atlanta, GA}
\email{rmickens@cau.edu}

\author[Washington]{Talitha Washington}
%\date{\today}
\address{Department of Mathematics, Howard University, Washington, DC}
\email{talitha.washington@howard.edu}

%\author[Initial Surname]{Author}
%\date{\today}
%\address{Address}
%\email{example@mail.com}

\maketitle

\let\thefootnote\relax
\footnotetext{AMS CLASSIFICATION: 35B09; 35G16; 35K05; 65M06; 80A10;}

%\begin{amscode} %new?
%35B09; 35G16; 35K05; 65M06; 80A10;
%\end{amscode}
%AMS 35B09, 35G16, 35K05; 65M06; 80A10; 8019

\begin{abstract}
We show that it is possible to construct microscopic-level discrete equations from macroscopic modeling PDEs for heat conduction in one space dimension. The significance of this result is that, in general, one starts from microscopic theories and then take their continuum limits to obtain the corresponding macroscopic PDEs, whereas here it is demonstrated that the reverse procedure is also possible. While our focus is on heat conduction, we discuss the applicability of our methodology to other physical systems.\\

\end{abstract}

\bigskip

\noindent {\Large \bf 1.~Introduction}\\

A physical system involving simple heat conduction may be examined from either macroscopic or microscopic perspectives \cite{ref3, ref5, ref14}. At the macroscopic level (MAL), the space-time variables are continuous, and the evolution equations are partial differential equations (PDEs) However, at the microscopic level (MIL), the space-time variables are generally discrete, and the governing equations are either difference equations (DE) or systems of coupled ordinary differential equations (ODE). For the heat transfer problem, a standard procedure is the assume the existence of a micro-level mathematical model and then take a ``continuum limit” to determine the macro-level PDEs \cite{ref3, ref10, ref12, ref14}. However, the reverse process, i.e., beginning with the MAL PDEs and using them to derive the MIL evolutionary DEs is generally thought to be either impossible or not able to provide a unique set of discrete MIL equations. A thorough search of the published research literature turns up no publications on this issue.

The major goal of this article is to demonstrate that from various MAL equations of simple heat transfer, it is possible to derive an essentially unique set of MIL equations. The existence of this result may provide insights as to whether this type of methodology can be applied to other important physical systems, such as wave propagation and fluid dynamics \cite{ref3, ref5, ref10, ref14}.

In outline, the paper proceeds as follows: Section 2 provides definitions and discussions related to both the MAL and MIL limit equations and how to determine one given the other. In Section 3, we start with three MAL PDEs and show how to select finite difference discretizations such that the MIL DEs correspond to random walk equations. Finally, Section 5 summarizes our results and discusses their implications for other physical systems.

It should be indicated that all of our current investigations, as they relate to this paper, refer to various versions of heat conduction for only one space dimension. Further, all parameters appearing in the associated PDEs are assumed to be constant. The dependent variable, $u$, is taken to represent the absolute Kelvin temperature, which has the property that it can not be negative. 

Throughout this paper, we will have the opportunity to use the following abbreviations:
\begin{itemize}
\item CL: continuum limit
\item DE: discrete equation
\item MAL: macroscopic level
\item MIL: microscopic level
\item MM: mathematical model
\item ODE: ordinary differential equation
\item PDE: partial differential equation\\
\end{itemize}

\noindent {\Large \bf 2.~Background}\\

For the purposes of this paper, a mathematical model of a physical system is an abstract representation of the system, which is expressed in terms of discrete equations (DE), ordinary differential equations (ODEs), or partial differential equations (PDEs). The dependent variables in these mathematical structures correspond to the dependent variables of the physical system, and it is assumed that the independent variables can be related to the usual space and time coordinates.

Microscopic or micro-level (MIL) mathematical models are generally formulated within a context in which an ``atomic substructure” of the physical system is assumed to exist. The mathematical model equations will be dependent on space, time, and system variables, which may be discrete and/or continuous-valued. Also, parameters may appear in these equations, and they will be assumed to be constant. (The parameters provide a characterization of the overall physical properties and behaviors of the MIL modeling system.) It should be noted that the MIL equations are usually DEs or systems of coupled ODEs \cite{ref5, ref10, ref14}.

Macroscopic or macro-level (MAL) equations for a physical system are usually derived by taking the continuum limit (CL) of the MIL equations \cite{ref5, ref10}. However, for a given set of MIL equations, this procedure may not yield unique mathematical expressions for the MAL equations, which in some cases are PDEs \cite{ref10, ref14}.

The CLs are obtained by letting some or all of the parameters in the MIL equations go to special limiting values, such that particular combinations of the parameters are held fixed \cite{ref5, ref10, ref14}. These combinations may usually be identified with measurable macroscopic physical properties of the system.

We now illustrate these concepts by considering the case of simple heat conduction in a one-space dimension \cite{ref14}. For this situation, a way to start is to assume that the MIC MM is a random walk equation

\begin{equation}
u_m^{k+1}=\left(\frac{1}{2}\right)u_{m+1}^k+\left(\frac{1}{2}\right)u_{m-1}^k \label{2pt1}
\end{equation}
where $(k, m)$ refers to discrete time and space variables that are associated with parameters $(\tau, a)$, which have the physical units, respectively, of time and distance. See Zaudere [\cite{ref14}, Chapter 1] for full explanatory details. 

If we make the definitions
\begin{subequations}
\begin{align}
k &\rightarrow t_k = \tau k, k = (0, 1, 2, …), \label{2pt2a} \\
m &\rightarrow x_m = am, m = \text{ integers}, \label{2pt2b}\\
u_m^k &\rightarrow u(x_m, t_k) \label{2pt2c}
\end{align}
\end{subequations}
then the continuum generalization of Eq.~(\ref{2pt1}) is
\begin{equation}
u(x,t+\tau)=\left(\frac{1}{2}\right)u(x+a, t)+\left(\frac{1}{2}\right)u(x-a,t). \label{2pt3}
\end{equation}
After some algebraic manipulations, this last equation can be rewritten to the form
\begin{equation}
\frac{u(x,t+\tau) - u(x,t)}{\tau} = \left(\frac{a^2}{2\tau}\right)\left[\frac{u(x+a,t)-2u(x,t)+u(x-a,t)}{a^2}\right]. \label{2pt4}
\end{equation}
If now we make the (really big) assumption that $u(x,t)$ has a Taylor series and take the limits 
\begin{equation}
a \rightarrow 0, \quad \tau \rightarrow 0, \text{ with } \frac{a^2}{2\tau} = D = \text{ fixed,} \label{2pt5}
\end{equation}
then from Eq.~(\ref{2pt4}), the following result is obtained
\begin{equation}
u_t=Du_{xx}, \quad \partial_t = \partial / \partial t, \quad \partial_{xx} = \partial^2/\partial x^2 \label{2pt6}
\end{equation}
The constant $D$ has the physical units of 
\begin{equation}
[D] = \frac{\text{Length}^2}{\text{Time}}. \label{2pt7}
\end{equation}
and is called the thermal diffusivity \cite{ref10, ref14}.

The just presented calculation is an example of how the macroscopic heat conduction PDE, Eq.~(\ref{2pt6}), can be ``derived” from the microscopic-based random walk equation. However, a number of other PDEs have been put forth to model heat conduction at the macroscopic level. Equations (\ref{3pt2}) and (\ref{3pt3}) are two such cases.

As stated previously, the central goal of this paper is to show that the inverse problem is not meaningless, i.e., the answer to the question, ``Given a MAL PDE modeling heat conduction, can we construct or derive a suitable MIL mathematical representation?”, is yes. Our general procedure is as follows:
\begin{itemize}
\item[(1)] First, select a MAL PDE.
\item[(2)] Using finite differences, construct a discretization of the PDE. This is a critical step and will require both insight and innovation into the various possibilities for discretizing derivatives.
\item[(3)] Rewrite the scheme as a linear combination of terms involving $u_m^{k+1}, u_{m+1}^{k+1}, u_{m-1}^{k+1}, u_m^{k-1}$, etc. The coefficients of these terms will depend on the parameters appearing in the PDE and the space and time step-sizes, $\Delta x$ and $\Delta t$. 
\item[(4)] Except for the terms involving $u_m^{k+1}, u_{m+1}^k, u_m^k, u_{m-1}^k$, set all the other coefficients to zero and determine if a mathematical consistent set of solutions exist for $\Delta x$ and $\Delta t$, expressed as functions of the PDE parameters.
\item[(5a)] If such a consistent set of solutions for $\Delta x$ and $\Delta t$ does not exist, then construct another discretization for the MAL PDE.
\item[(5b)] if solutions for $\Delta x$ and $\Delta t$ exist, then rewrite the resulting finite difference scheme in the form
\begin{equation}
u_m^{k+1} = p\left(u_{m+1}^k+u_{m-1}^k\right)+(1-2p)u_m^k \label{2pt8}
\end{equation}
and require
\begin{equation}
0<p\leq \frac{1}{2} \label{2pt9}
\end{equation}
note that in general, $p$ will depend on both $\Delta x$ and $\Delta t$, along with the PDE parameters, and provides a relationship between $\Delta x$ and $\Delta t$, once a value for $p$ is chosen.\\
\end{itemize}

\noindent {\large \bf 2.1~Comments:}\\

\begin{itemize}
\item[(i)] Note that Eq.~(\ref{2pt8}) is a discrete random walk equation \cite{ref14}.
\item[(ii)] The condition presented in Eq.~(\ref{2pt9}) ensures that $u_m^k$ satisfies a positivity restriction, i.e.,
\begin{equation}
u_m^k \geq 0 \Rightarrow u_m^{k+1} \geq 0, \quad m = \text{ integers.} \label{2pt10}
\end{equation}
\item[(iii)] Observe that the (iterative obtained) solutions to Eq.~(\ref{2pt8}) only require a knowledge of $u(x^m,0) = u_m^0$, $m = \text{integers}$, and the specification of the boundary conditions. This implies that we do not need to know $u_t(x,0)$, as required by some of the MAL PDE formulations of heat conduction, such as the Maxwell-Cattaneo equation \cite{ref1, ref8, ref9}.
\end{itemize}

In the next section, we will illustrate this methodology by applying it to three different formulations of the MAL PDEs for heat conduction. \\

\noindent {\Large \bf 3.~Three Heat Conduction PDEs}\\

Simple heat conduction in one-space dimension has been modeled at the macroscopic level by many sets of PDEs \cite{ref1, ref2, ref3, ref8, ref11, ref12, ref13, ref14}. The just listed references provide a detailed statement of the historical important articles on this topic; in particular, see the compilation provided by Sobolev \cite{ref12, ref13}.

In the calculations to be performed, the following three PDEs will be investigated.
\begin{subequations}
\begin{align}
u_t &= Du_{xx} \quad (\textrm{standard heat PDE \cite{ref14}}), \label{3pt1} \\
\tau u_{tt} + u_t &= Du_{xx} \quad (\textrm{Maxwell-Cattaneo model \cite{ref11, ref8}}), \label{3pt2}\\
u_t &= Du_{xx} +(\tau D)u_{txx} \quad (\textrm{symmetry derived model \cite{ref9}}). \label{3pt3}
\end{align}
\end{subequations}

It should be noted that $D$ is a constant diffusion coefficient, while $\tau$ is a parameter having the physical dimensions of ``time.” Also, the $\tau$'s appearing in Eqs.~(\ref{3pt2}) and (\ref{3pt3}) are not required to have the same numerical values.\\

\noindent {\Large \bf 4.~ Microscopic Limit Equations}\\

In this section, we will derive the microscopic or micro-level discrete equations for the three macroscopic or macro-level PDEs listed in Section 3.\\

\noindent {\large \bf 4.1~Standard Heat PDE}\\

A finite difference discretization of Eq.~(\ref{3pt1}) is \cite{ref14}
\begin{equation}
\frac{u_m^{k+1}-u_m^k}{\Delta t} = \left[\frac{u_{m+1}^k-2u_m^k+u_{m-1}^k}{(\Delta x)^2}\right], \label{4pt1}
\end{equation}
where the first-order time derivative and the second-order space derivatives have been replaced, respectively, by forward-Euler and central-difference schemes \cite{ref14}, and $(\Delta x, \Delta t)$ are the space and time step-sizes. The following notation is used
\begin{subequations}
\begin{align}
x &\rightarrow x_m = (\Delta x)m, \quad m = (\text{integers}) \label{4pt2a} \\
t &\rightarrow t_k = (\Delta t)k, \quad k=(0,1,2, ..) \label{4pt2b} \\
u(x,t) &\rightarrow u(x_m,t_k) = u_m^k \label{4pt2c}
\end{align}
\end{subequations}

After some minor algebraic manipulations, Eq.~(\ref{4pt1}) can be rewritten to the form
\begin{equation}
u_m^{k+1} = (1-2p)u_m^k + p(u_{m+1}^k + u_{m-1}^k), \label{4pt3}
\end{equation}
where
\begin{equation}
p = \frac{(\Delta t) D}{(\Delta x)^2}. \label{4pt4}
\end{equation}

If $p$ is restricted to the interval 
\begin{equation}
0<p \leq \frac{1}{2} \label{4pt4}
\end{equation}
then Eq.~(\ref{4pt3}) is a modified random walk equation. For $p=\frac{1}{2}$, we obtain the usual random walk discrete equation \cite{ref14}
\begin{equation}
u_m^{k+1} = \left(\frac{1}{2}\right)u_{m+1}^k+\left(\frac{1}{2}\right) u_{m-1}^k \label{4pt5}
\end{equation}
which is the starting point for deriving the standard heat PDE \cite{ref9, ref10, ref14} by taking the so-called continuum limit
\begin{equation}
\Delta x \rightarrow 0, \quad \Delta t \rightarrow 0, \quad \frac{(\Delta x)^2}{2(\Delta t)} = D = \text{ fixed}. \label{4pt6}
\end{equation}

The above analysis demonstrates that given the standard heat PDE, we can obtain a microscopic level discrete equation which is exactly the same as the random walk equation used to derive the PDE.

Observe that for this case, there is a functional relationship between the step-sizes, $\Delta x$ and $\Delta t$.\\

\noindent {\large \bf 4.2~Maxwell-Cattaneo Model}\\

This model PDE for heat conduction \cite{ref1, ref8} is also known as the damped wave equation \cite{ref14}. We select as our discretization the form 
\begin{equation}
\begin{split}
&\tau \left[\frac{u_m^{k+1}-2u_m^k+u_m^{k-1}}{(\Delta t)^2}\right] + \left[\frac{u_m^k-u_m^{k-1}}{\Delta t}\right] \\
&= D\left[\frac{u_{m+1}^k-u_m^{k+1}-u_m^{k-1}+u_{m-1}^k}{(\Delta x)^2}\right], \label{4pt7}
\end{split}
\end{equation}
where the first and second time derivatives are replaced, respectively, by central and backward-Euler discretizations and the second-order space derivative has a DuFort-Frankel discrete representation \cite{ref8, ref14}.Collecting together the various terms of the $(u_m^k)$ functions gives the expressions
\begin{align}
\begin{split}
\left[1+\left(\frac{D}{\tau}\right)\left(\frac{\Delta t}{\Delta x}\right)^2\right]u_m^{k+1} &= \left[2 - \frac{\Delta t}{\tau}\right]u_m^k
 + \left(\frac{D}{\tau}\right)\left(\frac{\Delta t}{\Delta x}\right)^2\left(u_{m+1}^k+u_{m-1}^k\right) \cr
&+\left[-1 + \frac{\Delta t}{\tau}-\left(\frac{D}{\tau}\right)\left(\frac{\Delta t}{\Delta x}\right)^2\right]u_m^{k-1}. \label{4pt8}
\end{split}
\end{align}

If we set
\begin{equation}
\Delta t = 2\tau, \label{4pt9}
\end{equation}
then the coefficient of the $u_m^k$ term is zero. Likewise, substituting this result into the coefficient of the $u_m^{k-1}$ term shows that it can be made zero if $\Delta x$ is selected to be
\begin{equation}
(\Delta x)^2 = 4D\tau. \label{4pt10}
\end{equation}
From these last two results, it follows that
\begin{equation}
1+\left(\frac{D}{\tau}\right)\left(\frac{\Delta t}{\Delta x}\right)^2 = 2. \label{4pt11}
\end{equation}
Consequently, Eq.~(\ref{4pt8}) can be finally written as
\begin{equation}
u_m^{k+1}=\left(\frac{1}{2}\right)u_{m+1}^k+\left(\frac{1}{2}\right)u_{m-1}^k. \label{4pt12}
\end{equation}
provided $\Delta x$ and $\Delta t$ satisfy the conditions given in Eqs.~(\ref{4pt9}) and (\ref{4pt10}). This is just the random walk equation. Note the important fact that from Eqs.~(\ref{4pt9}) and (\ref{4pt10}) that
\begin{equation}
\frac{(\Delta x)^2}{2\Delta t} = D, \label{4pt13}
\end{equation}
but no limits are to be taken, since both $\Delta x$ and $\Delta t$ are now determined by the two parameters, $\tau$ and $D$, occurring in the macroscopic level PDE given by Eq.~(\ref{3pt2}).

If we take $(\Delta x, \tau)$ to be micro-level parameters, and $(\tau, D)$ to be the corresponding macro-level parameters, then Eqs.~(\ref{4pt9}), (\ref{4pt10}), and (\ref{4pt13}) provide relations between their values. Finally, an easy calculation allows the determination of a length scale associated with the parameters, $\tau$ and $D$, it is 
\begin{equation}
L = \sqrt{\tau D} \label{4pt14}
\end{equation}
Observe that $\tau$ is both a microscopic and macroscopic variable. \\

\noindent {\large \bf 4.3~Another Standard Model PDE Discretization}\\

An interesting discretization of the standard heat PDE, $u_t = Du_{xx}$ is
\begin{equation}
\frac{u_m^{k+1}-u_m^{k-1}}{2 \Delta t} = D\left[\frac{u_{m+1}^k-u_m^{k+1}-u_m^{k-1}+u_{m-1}^k}{(Delta x)^2}\right], \label{4pt15}
\end{equation}
where the first-order time derivative is replaced by a central difference scheme and the second-order space derivative has the DuFort-Frankel discretization. This expression can be rewritten to the form
\begin{equation}
\left[1+\frac{2(\Delta t)D}{(\Delta x)^2}\right]u_m^{k+1}=\left[\frac{2(\Delta t)D}{(\Delta x)^2}\right](u_{m+1}^k+u_{m-1}^k)+\left[1-\frac{2(\Delta t)D}{(\Delta x)^2}\right]u_m^{k-1} \label{4pt16}
\end{equation}

If we select
\begin{equation}
\frac{2(\Delta t)D}{(\Delta x)^2}=1, \label{4pt17}
\end{equation}
then Eq.~(\ref{4pt16}) becomes
\begin{equation}
u_m^{k+1}=\left(\frac{1}{2}\right)u_{m+1}^k+\left(\frac{1}{2}\right)u_{m-1}^k, \label{4pt18}
\end{equation}
which is the random walk equation.\\

\noindent {\large \bf 4.4~Symmetry Derived Heat PDE}\\

Based purely on symmetry arguments, a heat conduction PDE can be derived \cite{ref9}; it is
\begin{equation}
u_t=Du_{xx}+(\tau D)u_{txx}, \label{4pt19}
\end{equation}
where the parameter, $\tau$, has the physical units of time.

The discretization is attained as follows:
\begin{itemize}
\item[(a)] Replace $u_t$ by the nonstandard representation
\begin{equation}
u_t \rightarrow \left[\left(\frac{u_{m+1}^{k+1}+u_m^{k+1}+u_{m-1}^{k+1}}{3}\right)-u_m^k\right] / \Delta t; \label{4pt20}
\end{equation}
\item[(b)] For $u_{xx}$ use
\begin{equation}
u_{xx} \rightarrow \frac{u_{m+1}^k-2u_m^k+u_{m-1}^k}{(\Delta x)^2}; \label{eq4pt21}
\end{equation}
\item[(c)] In $u_{txx}$, use a forward-Euler representation for the first-order time derivative and a central difference scheme for the second-order space derivative and obtain
\begin{equation}
u_{txx} \rightarrow \left[\frac{\tau D}{(\Delta t)(\Delta x)^2}\right](u_{m+1}^{k+1}-2u_m^{k+1}+u_{m-1}^{k+1}-u_{m+1}^k+2u_m^k-2u_{m-1}^k). \label{4pt22}
\end{equation}
\end{itemize}

Substituting these replacements into Eq.~(\ref{4pt19}) and collecting like terms in the $(u_m^k)$ together, gives the result
\begin{equation}
\begin{split}
&\left[1-\frac{3\tau D}{(\Delta x)^2}\right]\left(u_{m+1}^{k+1}+u_{m-1}^{k+1}\right)+\left[1+\frac{6\tau D}{(\Delta x)^2}\right]u_m^{k+1} = \cr
&\left[\frac{3(\Delta t)D}{(\Delta x)^2}-\frac{3\tau D}{(\Delta x)^2}\right]\left(u_{m+1}^k+u_{m-1}^k\right)+\left[3-\frac{6(\Delta t)D}{(\Delta x)^2}+\frac{6\tau D}{(\Delta x)^2}\right]u_m^k. \label{4pt23}
\end{split}
\end{equation}

If we require
\begin{equation}
\frac{3\tau D}{(\Delta x)^2} = 1, \quad \Delta t=2\tau, \label{4pt24}
\end{equation}
then Eq.~(\ref{4pt23}) becomes
\begin{equation}
u_m^{k+1}=\left(\frac{1}{3}\right)u_{m+1}^k+\left(\frac{1}{3}\right)u_m^k+\left(\frac{1}{3}\right)u_{m-1}^k, \label{4pt25}
\end{equation} 
and this is a form of the random walk equation \cite{ref14}. Also, note that
\begin{equation}
\Delta x=\sqrt{3 \tau D}. \label{4pt26}
\end{equation}

\

\noindent {\Large \bf 5.~Discussion}\\

In general, the investigation of the properties of a physical system begins with the formulation of a mathematical model centered on its MIL structure. However, few systems have MMs for which good analytical approximations or numerical solutions exist for the MIL equations. A partial solution to resolve this issue is to examine special limiting forms of the MIL equations that allow analytical and/or numerical solutions to be determined. These derived MAL equations are usually PDEs, in contrast to the original MIL equations which are either DEs or coupled ODEs \cite{ref3, ref5, ref9, ref10}. The MAL equations are constructed by taking CLs of the MIL equations \cite{ref10, ref14}.

An issue of concern is that while going from the MIL equations to the MAL equation has been studied extensively, little or no effort has been made to go in the opposite direction, i.e., derive MIL equations from MAL equations. In fact, our search of the research literature has found no such work. The major goal of this paper is to show that this task is possible, at least when it concerns simple heat transfer as modeled by the standard heat PDE. Another way of stating this problem is to ask the question: Can we derive general microscopic theories (MIL equations) from just information provided by the MAL mathematical equations? To repeat, our answer is yes and the calculations given in Section 5 demonstrate this to be correct.

For simple heat conduction, we conclude that the various PDE models at the MAL all contain remnants of an underlying classical physics ``atomonistic structure" and this structure can be mathematically reconstituted by the process of discretization with the imposition of several other constraints.

Furthermore, our work provides an illustration of the fact that discretization of PDEs may be performed for reasons unrelated to the construction of schemes to calculate numerical solutions. In fact, it is often the case that the DEs are better mathematical representations of physical systems than the corresponding differential equations \cite{ref6, ref7, ref12, ref13, ref14}. This result is related to the realization that all experiments are carried out under conditions where the space and time variables are measured at discrete values, and the values of dependent variables are always approximated by rational rather than real numbers \cite{ref3, ref5, ref10}.

It should be indicated that our general methodology can be applied to determine the MIL equations for the wave equation
\begin{equation}
u_{tt}=c^2u_{xx}, \quad u= u(x,t), \label{5pt1}
\end{equation}
where $u(x,t)$ is the transverse displacement of a one-space dimension string at location $x$ and time, $t$. This may be shown by defining $u_m(t)$ as
\begin{equation}
u_m(t)=u(x_m, t), \quad x_m=(\Delta x)m, \quad m = (\text{integers}). \label{5pt2}
\end{equation}

Consequently, discretizing the right-side of Eq.~(\ref{5pt1}) gives
\begin{equation}
\frac{d^2u_m}{dt^2}=\left[\frac{c^2}{(\Delta x)^2}\right](u_{m+1}-2u_m+u_{m-1}). \label{5pt3}
\end{equation}
Note that this expression can be rewritten to the form
\begin{equation}
M\frac{d^2u_m}{dt^2}=k(u_{m+1}-u_m)-(u_m-u_{m-1}), \label{5pt4}
\end{equation}
where
\begin{equation}
k=\frac{Mc^2}{(\Delta x)^2}, \label{5pt5}
\end{equation}
and $M$ is the mass of a particle located at $x_m$. One interpretation of Eq.~(\ref{5pt4}) is that it represents a chain of identical coupled simple linear harmonic oscillators having for each oscillator the effective mass, $M$, and force constant, $k$. This is one of the standard models used to derive the wave PDE starting at the MIL with Eq.~(\ref{5pt4}) \cite{ref10, ref14}.

\begin{figure}
\centering
\includegraphics[scale=0.45]{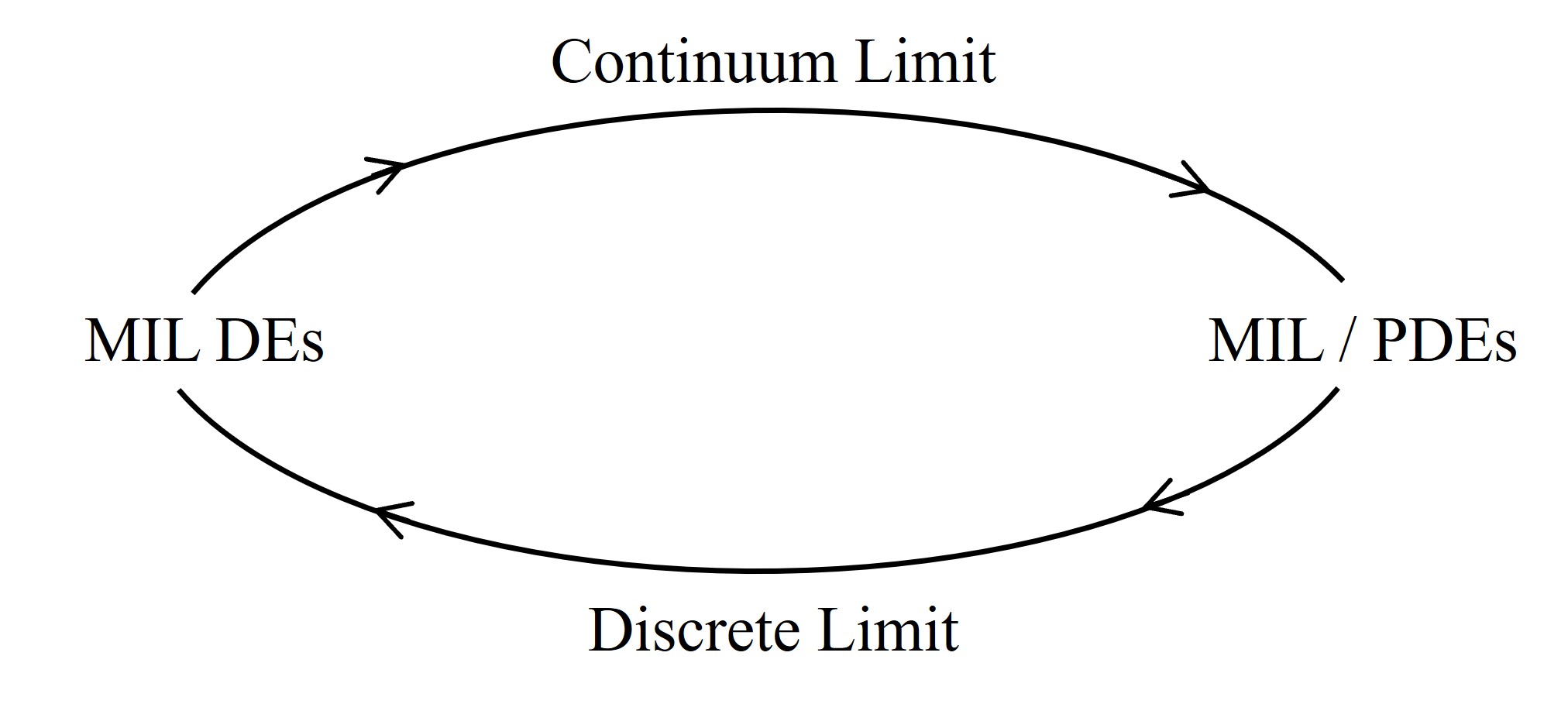}
\caption{Connections between the micro-level (MIL) and macro-level (MAL) mathematical models for a physical system. The corresponding equations are, respectively, discrete (DE) or Partial differential equations (PDEs). The arrow indicate the methodologies needed to transition between MIL and MAL equations.}
 \label{fig1}
\end{figure}

It would be of immense value to the physical and engineering sciences if the methodology of this paper could be extended to MAL nonlinear PDEs. If this can be accomplished, then some of the issues raised in this classic paper by Lagerstam \cite{ref4} might be resolved. See Figure \ref{fig1} for a pictorial illustration of what we have discussed in this paper. It shows that, contrary to what is generally believed, the CL PDEs may be used to reconstruct the underlying DEs for the associated MIL equations. For heat conduction problems, the DEs will be mathematical expressions which can be interpreted as modified random walk equations.

Finally, we end this paper with the following comments that should be be forgotten by researchers who investigate issues involving heat conduction. The MIL models are not, in general, representations of the actual physical system. They only represent MMs that reproduce certain desired physical properties of the phenomena of interest to our investigations. For example, no real heat transfer issue is an actual random walk of ``atoms” at the MIL, where the ``particles" move along a one-dim space lattice. But taking such a MM model provides a MAL PDE (through taking the CL) where its relevant solutions predict results in general agreement with experimental results.\\

\noindent {\Large \bf 6.~Acknowledgements}\\

The authors thank Imani Beverly and Brian Briones, reference librarians at the Atlanta University Center Robert W.~Woodruff Library for their help in obtaining documents, articles, and books related to this research.\\

Author Contributions: This manuscript contains equal contributions from both authors.\\

Financial Disclosure: The work of this manuscript was not supported by external funding sources.\\

Conflict of Interest: The authors declare no potential, nor actual conflict(s) of interest.\\

%\noindent {\Large \bf 8.~References}

%\begin{thebibliography}{99}

%\bibitem{1} Spiegel, M. R. (1981). Theory and problems of Advanced Calculus: Si (metric) edition. McGraw-Hill. 

%\end{thebibliography}
\end{document}